\documentstyle[12pt,titlepage,twoside]{article}
\font\twelvegtc=eufm10 scaled 1200

\font\ninegtc=eufm9
\font\sevengtc=eufm7

\newfam\gtcfam

\textfont\gtcfam=\twelvegtc
\scriptfont\gtcfam=\ninegtc
\scriptscriptfont\gtcfam=\sevengtc
\font\twelveBBB=msbm10 scaled 1200
\font\tenBBB=msbm10
\font\sevenBBB=msbm7
\newfam\BBBfam
\def\BBB{\fam\BBBfam\twelveBBB}
\textfont\BBBfam=\twelveBBB
\scriptfont\BBBfam=\tenBBB
\scriptscriptfont\BBBfam=\sevenBBB
\def\QQ{{\BBB Q}}   \def\ZZ{{\BBB Z}}    
   \def\CC{{\BBB C}}    \def\PP{{\BBB P}}

\def\VEC#1,#2{(#1_1,#1_2,\dots,#1_{#2})}
\def\OVEC#1,#2{(#1_0,#1_1,#1_2,\dots,#1_{#2})}
\def\SET#1,#2{\{#1_1,#1_2,\dots,#1_{#2}\}}
\def\FAM#1,#2{\ #1_1,#1_2,\dots,#1_{#2}\ }
\def\BSER#1,#2,#3{#1_0+#1_1#2+\cdots+#1_{#3}#2^{#3}}
\def\SER#1,#2{#1_0+#1_1#2+#1_2#2^2+\cdots}
\def\POL#1,#2,#3{#1_0#2^{#3}+#1_1#2^{#3-1}+\cdots+#1_{#3-1}#2+#1_{#3}}
\def\UPOL#1,#2,#3{#2^{#3}+#1_1#2^{#3-1}+\cdots+#1_{#3-1}#2+#1_{#3}}

\def\rp#1,#2{{\rm Rp}_{#1}#2}
\def\lp#1,#2{{\rm Lp}_{#1}#2}
\def\bl{\!\in\!}

\def\hom{{\rm Hom}}
\def\ext{{\rm Ext}}
\def\sgn{{\rm sgn}\,}
\def\fa#1{\;\forall#1\;}
\def\fain#1,#2{\fa#1\bl#2}
\def\map#1,#2{#1\longrightarrow#2}
\def\MAP#1,#2,#3{#1\,\colon\;\map#2,#3}
\def\sp#1,#2{\langle#1,#2\rangle}
\def\comsq#1,#2,#3,#4,#5,#6,#7,#8{
\begin{array}{ccc}
     #1&\stackrel{#2}{\longrightarrow}&#3\\
     \llap{$\scriptstyle #4$}\Biggl\downarrow\Biggr.& &
                      \Biggl\downarrow\Biggr.\rlap{$\scriptstyle #5$}\\
     #6&\stackrel{#7}{\longrightarrow}&#8
\end{array}}
\def\ldf#1,#2{
\par\noindent\hbox to\textwidth
{#1\hfil #2\hfil}\par
}
\def\rdf#1,#2{
\par\noindent\hbox to\textwidth
{\hfil #2\hfil #1}\par
}
\newcounter{No}   \newcounter{SubNo}[No]    \newcounter{SubSubNo}[SubNo]
\renewcommand{\theNo}{\S\arabic{No}}
\renewcommand{\theSubNo}{\arabic{No}.\arabic{SubNo}}
\renewcommand{\theSubSubNo}{\arabic{No}.\arabic{SubNo}.\arabic{SubSubNo}}

\def\No#1{\refstepcounter{No}\par\vfill\pagebreak[3]\noindent
                               {\large\bf\theNo.\hspace{2pt}#1.}\par}
\def\SubNopt#1{\refstepcounter{SubNo}\vspace{2ex}\par\noindent
                               {\bf\theSubNo.\hspace{2pt}#1.}\hspace{1ex}}
\def\SubNo#1{\refstepcounter{SubNo}\vspace{2ex}\par\noindent
                               {\bf\theSubNo.\hspace{2pt}#1}\hspace{1ex}}
\def\SubSubNo#1{\refstepcounter{SubSubNo}\vspace{1ex}\par\noindent
                               {\bf\theSubSubNo.}\hspace{2pt}#1.\hspace{1ex}}
\def\EF{\endgroup\par\vspace{1ex}\par}

\def\Pr{\SubSubNo{PROPOSITION}\begingroup\sl}

\def\Cl{\SubSubNo{COROLLARY}\begingroup\sl}
\def\proof{\par\noindent{\sc Proof.}\hspace{1ex}}
\def\Cs{\begin{enumerate}}
\def\EC{\end{enumerate}}
\def\Th{\SubSubNo{THEOREM}\begingroup\sl}
\def\Lm{\SubSubNo{LEMMA}\begingroup\sl}

\title{Semistable sheaves on Del Pezzo surfaces and Kronecker modules}
\author{B.V. Karpov \thanks{Preliminary work on this paper in Moscow was
supported by
the foundation PRO MATHEMATICA (France) and the J.Sorros foundation (USA)}\\
Algebra and Math.Logics Department\\ of the Moscow Institute of Electronics and
Mathematics\\
109028 Trehsvjatitel'skij (ranee B.Vuzovskij)\\ per. 3/12, Moscow, Russia \\
e-mail: root@aml.miemstu.msk.su}
\date{September 1994}
\begin{document}
\maketitle
\begin{abstract}
In this paper we deal with semistable sheaves which can be represented as the
cokernel of an injective (or as the kernel of a surjective) morphism
$E_1\otimes\protect\CC^m\longrightarrow E_2\otimes\protect\CC^n$ , where
$E_1$ and $E_2$ are exceptional bundles. To each such a sheaf we assign the
linear map $\protect\CC^m\otimes\hom^*(E_1,E_2)\longrightarrow
\protect\CC^n$. We obtain sufficient conditions on the topological
invariants of the sheaves for the moduli space of the sheaves to be
isomorphic to the quotient of an open subset in
$\protect\PP{\cal L}(\protect\CC^m\otimes\hom^*(E_1,E_2),\protect\CC^n)$
under the action of SL$(\protect\CC^m)\times$SL$(\protect\CC^n)$.
\end{abstract}
%
%
\markboth{B.V.Karpov}{Semistable sheaves on Del Pezzo surfaces}
\renewcommand{\theenumi}{(\roman{enumi})}
%
%
\No{Introduction}
For the first time exceptional vector bundles appeared in the work [DrLP]
of Drezet and Le Potier in 1985. The authors described the set of all triples
$(r,c_1,c_2)$ such that there exists a semistable sheaf $V$ on $\CC\PP_2$ with
the given rank and Chern classes:
          $$r(V)=r;\quad c_1(V)=c_1;\quad c_2(V)=c_2.
$$
Exceptional bundles appeared in this work as stable bundles with the
discriminant
$\Delta < \frac{1}{2}$. The ranks and Chern classes of exceptional bundles were
used for constructing the bound for these invariants of semistable sheaves.

It was proved that the bundle $E$ is exceptional if and
only if it is stable and rigid ($\ext^1(E,E)=0$). So any exceptional bundle $E$
has the moduli space consisting of one point.

Further on the concept of exceptional bundle was developed in two directions.
 From one hand, Gorodentsev introduced exceptional collections of bundles and
helices
 on $\PP_n$ (see [Go1]). It was founded that for any full exceptional
collection
there exists a spectral sequence of Beilinson type. Besides, the braid group
acts
on the set of all helices and on the set of all exceptional collections (a
survey of the helix theory and further references see in [Go3]).

 From the other hand, Drezet used exceptional bundles for obtaining more
concrete
information about moduli spaces of semistable sheaves on $\PP_2$ (see
[Dr1], [Dr2], [Dr3]).
In the work [Dr2] so-called semistable sheaves of height $0$ were considered.
The definition of height is specific for $\PP_2$ but the height of a semistable
sheaf $V$ equals $0$ if and only if there exists one of the following exact
triples:
$$0\longrightarrow E_1\otimes\protect\CC^m\longrightarrow
E_2\otimes\protect\CC^n
\longrightarrow V\longrightarrow 0
$$
or
$$0\longrightarrow V\longrightarrow E_1\otimes\protect\CC^m\longrightarrow
E_2\otimes\protect\CC^n\longrightarrow 0
$$
Here $E_1$ and $E_2$ are exceptional bundles on $\PP_2$. Suppose
$$\MAP t,\CC^m\otimes\hom^*{(E_1,E_2)},\CC^n
$$
is the linear map naturally associated with the morphism of bundles
$E_1\otimes\CC^m\longrightarrow E_2\otimes\CC^n$ in one of the exact triples;
then $t$ corresponds to $V$. Using this correspondence, Drezet has shown that
the moduli space of semistable sheaves of height $0$ is isomorphic to the
quotient of the open subset in $\PP{\cal
L}(\CC^m\otimes\hom^*(E_1,E_2),\CC^n)$,
consisting of semistable points with respect to the action of
SL$(\CC^m)\times$SL$(\CC^n)$, under this action.

Such moduli spaces of semistable sheaves on $\PP^1\times\PP^1$ were considered
in the work [Ka].
The subject of this paper is the results of such type for semistable sheaves on
Del Pezzo surfaces (see theorems 2.7.1 and 2.7.2).

I am grateful to the Max-Planck-Institut f\"ur Mathematik in Bonn for the
stimulating atmosphere in which this work was finished.

\No{General facts and formulation of main results}
\SubNopt{Notation}
Let $S$ be Del Pezzo Surface, $K$ be its canonical class. For any coherent
sheaf $E$ over $S$ with $r(E)>0$ put by definition
$$\nu(E)=\frac{c_1(E)}{r(E)}\in{\rm Pic}(S)\otimes\QQ$$
$$\mu (E)=\frac{c_1(E)\cdot (-K)}{r(E)}\in\QQ\mbox{ --- the {\it slope\/} of }
E
$$
$$q(E)=\frac{ch_2(E)}{2r(E)}=\frac{1}{2r(E)}(c^2_1(E)-2c_2(E))$$

For any two coherent sheaves $E$ and $F$ denote
$$\chi(E,F)=\sum (-1)^i\dim\ext^i(E,F).$$

The Hirzebruch -- Riemann -- Roch theorem gives the following expression
of $\chi$ in terms of invariants introduced above.
$$\chi(E,F)=r(E)r(F)\left( 1+\frac{1}{2}(\mu(F)-\mu(E))+q(F)+q(E)-
             \nu(E)\cdot\nu(F)\right)
$$
\noindent Sometimes we call this expression {\it Riemann -- Roch formula}.
One can easily see that after disclosing the parentheses the right hand side
is well defined as $r(E)r(F)=0$.

It is convenient for us to consider the {\it Mukai lattice\/}
$${\rm Mu}(S) = \ZZ\oplus{\rm Pic}(S)\oplus\ZZ.
$$
Any element $v=(r,c,s)\in{\rm Mu}(S)$ is called {\it Mukai vector\/}.
To each coherent sheaf $E$ on $S$ we assign Mukai vector
$$v(E)=(r(E),c_1(E),ch_2(E))$$

The values $\mu$ and $q$ are regarded as functions Mu$(S)\setminus\{ r=0\}
\longrightarrow\QQ$,
and $\nu$ is considered as a map Mu$(S)\setminus\{ r=0\}\longrightarrow{\rm
Pic}(S)\otimes\QQ$.
The Riemann -- Roch formula gives the bilinear form $\chi(*,*)$ on Mu$(S)$.
Sometimes we shall write $\chi(E,v)$ instead of $\chi(v(E),v)$.

We shall need also the following expression for the skew part $\chi_-$ of
the form $\chi$:
$$\begin{array}{ccccc}\vspace{2ex}\chi_- (E,F)&=&\chi(E,F)\ -\ \chi(F,E)&=& \\
                          &=&\left| \begin{array}{cc} \vspace{1ex}
                                   c_1(F)\cdot (-K)& r(F) \\
                                   c_1(E)\cdot (-K)& r(E)
\end{array} \right| &=& r(E)r(F)(\mu(F)-\mu(E))
\end{array}
$$
\SubSubNo{SWING LEMMA}\begingroup\sl
Suppose $0\longrightarrow E\longrightarrow F
\longrightarrow G\longrightarrow 0$
is exact triple of coherent sheaves, then:\\
a) $\chi_-(E,F)=\chi_-(F,G)=\chi_-(E,G)$\\
b) if $r(E)r(G)>0$ then there is one of the possibilities:
\begin{eqnarray*} &\mu(E)<\mu(G)<\mu(F)&\\
    {\sl or}\quad &\mu(E)=\mu(G)=\mu(F)&\\
    {\sl or}\quad &\mu(E)>\mu(G)>\mu(F)&.
\end{eqnarray*}
\EF
\noindent The proof is by direct calculation.
\SubNopt{Exceptional sheaves and exceptional collections}
A sheaf $E$ is {\it exceptional\/} if
$$\hom(E,E)=\CC,\qquad\ext^i(E,E)=0\quad{\rm for }\ i>0\ .$$

In the work [KuOr] it is proved that any exceptional sheaf on a Del Pezzo
surface
 is either locally free or a
 torsion sheaf of the form ${\cal O}_e(d)$, where $e$ is irreducible rational
curve with $e^2=-1$.

An ordered collection of sheaves $(E_1,\cdots ,E_n)$ is {\it exceptional} if
all $E_j$ are exceptional and $\ext^i(E_k,E_j)=0, \forall i ,$ for $k>j$.
Exceptional collection is {\it full} if it generates the bounded derived
category $D^bSh(S)$ of coherent sheaves on $S$.

\SubNopt{Stability}
We consider the usual Mumford -- Takemoto stability with
respect to the anticanonical class. Namely, a torsion free sheaf $E$ is said
to be {\it stable\/} ({\it semistable\/}) if for any subsheaf $F$ of $E$ such
that $r(F)<r(E)$ the following inequality holds:
$$\mu(F)<\mu(E)\qquad( \mu(F)\leq\mu(E)\quad\mbox{for semisatbility}  ).$$

Let us recall the following facts.
\Pr \\1. If $E$ and $F$ are semistable and $\mu(E)>\mu(F)$ then $\hom(E,F)=0$\\
2. Exceptional bundles are stable.
\EF

The first  assertion is the standard stability property, the second
one is proved in the work [KuOr].
\SubNopt{Mutations of exceptional pairs}
Let $(E,F)$ be an exceptional pair; then $\chi(F,E)=0$ and
$\chi(E,F)=\chi_-(E,F)$.
It is known (see [KuOr]) that the pair $(E,F)$ belongs to one
of the following types that are distinguished as follows:
\begin{itemize}
\item {\it hom} - pair: $\hom(E,F)\not=0, \ext^i(E,F)=0, i>0\Longleftrightarrow
\chi_-(E,F)>0$;
\item {\it ext} - pair: $\ext^1(E,F)\not=0, \ext^i(E,F)=0,
i=0,2\Longleftrightarrow
\chi_-(E,F)<0$;
\item {\it zero} - pair: $\ext^i(E,F)=0\quad$ for all
$i\quad\Longleftrightarrow\quad
\chi_-(E,F)=0$.
\end{itemize}

Let $(E,F)$ be a {\it hom} - pair. It is proved in [KuOr] that the canonical
map
$E\otimes\hom(E,F)\longrightarrow F$ is either surjective or injective. Hence,
there
exists one of the following exact triples:
\begin{equation}
0\longrightarrow L_E F \longrightarrow E\otimes\hom(E,F)\longrightarrow F
\longrightarrow 0
\label{lr}
\end{equation}
or
\begin{equation}
0\longrightarrow E\otimes\hom(E,F)\longrightarrow F\longrightarrow
L_E F\longrightarrow 0
\label{ls}\end{equation}
In the both cases the pair $(L_E F,E)$ is exceptional, it is called the {\it
left
 mutation\/} of the pair $(E,F)$. In the first case the left mutation is said
to be
{\it regular\/}, in the second one the left mutation is {\it singular\/}.

The right mutations of the pair $(E,F)$ are defined in the dual manner. Namely,
it is known that the partial dualization $E\longrightarrow F\otimes\hom(E,F)^*$
of
the canonical map is also either injective or surjective. Therefore, we have:
\begin{equation}
0\longrightarrow E\longrightarrow F\otimes\hom(E,F)^*\longrightarrow
R_F E\longrightarrow 0
\label{rr}\end{equation}
or
\begin{equation}
0\longrightarrow R_F E \longrightarrow E\longrightarrow F\otimes\hom(E,F)^*
\longrightarrow 0
\label{rs}\end{equation}
As above, the pair $(F,R_F E)$ is exceptional, it is called the {\it right
 mutation\/} of the pair $(E,F)$. In the first case the right mutation is
{\it regular\/}, in the second one it is {\it singular\/}.

Now let $(E,F)$ be an {\it ext} - pair. Left and right mutations of it are
the universal {\it extensions}:
\begin{equation}
0\longrightarrow F\longrightarrow L_E F\longrightarrow E\otimes\ext^1(E,F)
\longrightarrow 0
\label{le}\end{equation}
and
\begin{equation}
0\longrightarrow F\otimes\ext^1(E,F)^*\longrightarrow R_F E\longrightarrow E
\longrightarrow 0
\label{re}\end{equation}
As in the previous cases, the pairs $(L_E F,E)$ and $(F,R_F E)$ are
exceptional.

It can be easily seen that applying the definitions of mutations to a
{\it zero} - pair, we have always the permutation of its members.
\SubNopt{System of sheaves generated by a pair} Consider an exceptional pair
$(E,F)$
and an infinite set of sheaves $\{ E_i\}_{i\in\ZZ}\quad$ such that
$$E_1=E,\quad E_2=F\quad {\rm and}\quad E_{i+1}=R_{E_i}E_{i-1}\quad (\quad
               \Longleftrightarrow\quad  E_{i-1}=L_{E_i}E_{i+1})
$$
The set $\{ E_i\}$ is called the {\it system of sheaves generated by the pair}
$(E,F)$.

Applying the functor
$\hom(\ \cdot\ ,E)$ to the triples (\ref{lr}), (\ref{ls}), and (\ref{le}),
we see that the non-trivial space $\ext^{\textstyle{\cdot}}(E,F)$ is dual to
the non-trivial space $\ext^{\textstyle{\cdot}}(L_E F,E)$. Similarly, applying
the functor
$\hom(F,\ \cdot\ )$ to the triples (\ref{rr}), (\ref{rs}), and (\ref{re}),
we see that the non-trivial space $\ext^{\textstyle{\cdot}}(E,F)$ is dual to
the non-trivial space $\ext^{\textstyle{\cdot}}(F,R_F E)$. Therefore, we have
proved
the following
\Lm The non-zero dimension of the Ext-space between the neighbouring members of
the
system $\{ E_i\}$ (i.e. $\dim\hom(E_i,E_{i+1})$ if $(E_i,E_{i+1})$ is a
 {\it hom} - pair and
$\dim\ext^1(E_i,E_{i+1})$ if the type of the pair is {\it ext\/}) does not
depend
on $i$.
\EF
 We denote this dimension by $h$. It is obvious that
$$h=\vert\chi(E_i,E_{i+1})\vert\ .$$

If $h>1$, then there are the following possibilities for the system $\{ E_i\}$.
\begin{description}
\item{$(+)$}
All pairs $(E_i,E_{i+1})$ are of type {\it hom\/}. In this case all
$E_i$ are locally free.
\item{$(-)$}
There is only one {\it ext} - pair $(E_p,E_{p+1})$ , all the others are {\it
hom}
- pairs. In this case all $E_i,i\not=p$ are locally free. The sheaf $E_p$ can
be
either locally free or isomorphic to ${\cal O}_e(d)$.
\end{description}

This classification is proved in 3.1. Sometimes we write briefly
"$\{ E_i\}$ is of type $(+)$" or simply "$\{ E_i\}$ is $(-)$" meaning
 that the corresponding
possibility takes place for the system $\{ E_i\}$.

\vspace{1ex}

If $h=1$, then up to renumbering of $\{ E_i\}$ all the mutations are described
by the following exact triple
$$0\longrightarrow E_0\longrightarrow E_1\longrightarrow E_2
\longrightarrow 0$$
For instance, the left mutation of $(E_1,E_2)$ is regular, the right mutation
of
$(E_1,E_2)$ is singular, and $E_3\cong E_0$. Actually, $E_i\cong E_j$ iff
$i\equiv j$ mod 3. The simplest example of this situation is
$$0\longrightarrow {\cal O}(-e)\longrightarrow {\cal O}
\longrightarrow {\cal O}_e\longrightarrow 0\quad ,
$$
where $e$ is irreducible rational curve with $e^2=-1$.

\vspace{1ex}

In the case $h=0$ it is obvious that all pairs $(E_i,E_{i+1})$ have type {\it
zero},
$E_i\cong E_1$ whenever $i$ is even, and $E_i\cong E_2$ whenever $i$ is odd.

\vspace{1ex}

Suppose $h>2$; then there are two limits of the sequence $\{ \mu(E_i)\}$:
$$\mu_{-\infty}=\lim_{i \to {-\infty}} \mu(E_i)\quad {\rm and}\quad
\mu_{+\infty}=\lim_{i \to {+\infty}} \mu(E_i)
$$
(see proof of 3.1.4). If the system $\{ E_i\}$ has type $(+)$, then the slopes
of $E_i$ are ordered:
$$\mu_{-\infty}<\mu(E_i)<\mu(E_{i+1})<\mu_{-\infty}\quad \mbox{for all } i
$$
In the $(-)$ case the ordering is broked by the {\it ext} - pair
$(E_p,E_{p+1})$:
$$\mu(E_{p+1})<\mu(E_{p+2})<\cdots
<\mu_{+\infty}<\mu_{-\infty}<\cdots<\mu(E_{p-2})
<\mu(E_{p-1})\lbrack<\mu(E_p)\rbrack$$
(the last inequality is in brackets because $E_p$ can be of zero rank.)
\SubNopt{Kronecker moduli spaces}
Let $H_0,L,H_1$ be finite-dimensional vector spaces with $\dim L>2$ . Any
linear
map $\MAP {t},{H_0\otimes L},{H_1}$ is called a {\it Kronecker $L$-module\/}
(this term was used by Drezet in his work [Dr2]). A morphism from the
Kronecker module $t$ to the Kronecker module $\MAP {t'},{H'_0\otimes L},{H'_1}$
 is a pair of linear maps
$(f_0,f_1)$ , where $\MAP {f_0},{H_0},{H'_0}$ , $\MAP {f_1},{H_1},{H'_1}$ ,
such that the following diagram
$$\begin{array}{ccc}
H_0\otimes L&\stackrel{t}{\longrightarrow}&H_1\\
\llap{$\scriptstyle f_0\otimes {\rm id}_L$}\Biggl\downarrow\Biggr.& &
\Biggl\downarrow\Biggr.\rlap{$\scriptstyle f_1$}\\
H'_0\otimes L&\stackrel{t'}{\longrightarrow}&H'_1
\end{array}$$
is commutative.

It is clear that all the Kronecker $L$-modules form an Abelian category.
The Kronecker submodules, quotient-modules, and the kernels and cokernels
of morphisms are defined in the obvious way.

Let us denote $W={\cal L}(H_0\otimes L,H_1)\setminus\{ 0\}$. The reductive
group
$G_0=({\rm GL}(H_0)\times{\rm GL}(H_1))/ \CC^*$ acts on $W$ as follows:
$$(\CC^*(g_0,g_1),t)\longmapsto g_1\circ t\circ(g_0\otimes {\rm id}_L)^{-1}.
$$
Consider the projectivization $\PP=\PP({\cal L}(H_0\otimes L,H_1))$. The action
of $G_0$ on $W$ induces an action of the group ${\rm SL}(H_0)\times{\rm
SL}(H_1)$
on $\PP$.

A nonzero Kronecker module is {\it semistable\/} ({\it stable\/}) if its
image in $\PP$ is semistable (stable) in the usual sense of stability of a
point in an algebraic variety with reductive group action relative to
linearization (see [MuFo],[Dr2]). We shall use the following
criterion (see [Dr2]).
\Pr Let $\MAP {t},{H_0\otimes L},{H_1}$ be a nonzero Kronecker module.
Then $t$ is semistable (stable) if and only if for any submodule
$\MAP {t'},{H'_0\otimes L},{H'_1}$ of $t$ such that $H'_0\not=\{ 0\}$ and
$H'_1\not= H_1$ the following inequality holds:
$$\frac{\dim H'_1}{\dim H'_0}\ge\frac{\dim H_1}{\dim H_0}\quad\left(
{\it resp. }\frac{\dim H'_1}{\dim H'_0}>\frac{\dim H_1}{\dim H_0}\right)
$$
\EF

Denote by $W^{ss}$ ($W^s$) the set of all semistable (respectively stable)
Kronecker modules and put $h=\dim L$, $m=\dim H_0$, $n=\dim H_1$. Denote
also
$$N(h,m,n)=W^{ss}/G_0\quad {\rm and}\quad
  N^s(h,m,n)=W^{s}/G_0$$
These quotients are the {\it Kronecker moduli spaces}. From definitions it
follows that $N(h,1,1)=\PP_{h-1}$ and $N(h,1,n)=G(h-n,h)$; thus, the
Kronecker moduli spaces are in a sense generalizations of Grassmanians.
The first non-trivial case $N(3,2,2)\cong\PP^5$ is considered in the
work [Dr2].

Let $(E,F)$ be an exceptional pair. We denote by $ev_{E,F,x}$ the fiber
at a point $x\in S$ of the canonical map:
$$\MAP{ev_{E,F,x}},{E_x\otimes\hom{(E,F)}},{F_x}
$$
regarded as Kronecker $\hom(E,F)$-module. The following properties are proved
in
[Ka].

\Pr If the left mutation of  $(E,F)$ is regular, then
$$ev_{E,F,x}\mbox{ is (semi)stable}\quad\Longleftrightarrow\quad
     ev_{L_E F,E,x}\mbox{ is (semi)stable}$$
\EF
\Cl If the system $\ \{ E_i\}\ $ is of type $(+)$ and contains a
one-dimensional bundle,
then $ev_{E_i,E_{i+1},x}$ is stable for any $i\in\ZZ,\ x\in S$.
\EF
\SubNo{The main results} of this paper are the two following theorems.
\Th Suppose $(E_1,E_2,F_2,\cdots,F_l)$ is a full exceptional collection
of bundles over Del Pezzo surface $S$, the system $\{ E_i\}$ generated by the
pair $(E_1,E_2)$ has type $(-)$, and $h=\vert\chi(E_i,E_{i+1})\vert >2$.
Let $v=(r,c,s)$ be a Mukai vector with $r>0$ such that the following conditions
are
satisfied:

\vspace{1ex}
\ldf{$(0)$},{$\chi(F_j,v)=0,\quad j=2,\cdots,l$}
\vspace{1ex}
\ldf{$(1)$},{$\max\limits_j\mu(F_j)-K^2<\mu(v)<\min\limits_j\mu(F_j)$}
\ldf{$(2-)$},{$\mu_{+\infty}<\mu(v)<\mu_{-\infty}$}
\vspace{1ex}

Then for  $m=\vert\chi(E_3,v)\vert$ and
 $n=\vert\chi(E_2,v)\vert$
the Kronecker moduli space $N(h,m,n)${ }$(N^s(h,m,n))$ is isomorphic to
the coarse moduli space of semistable (resp. stable) sheaves with the
Mukai vector $v$.
\EF
\Th Let $(E_1,E_2,F_2,\cdots,F_l)$ be full exceptional collection
of bundles on $S$ such that the system $\{ E_i\}$ generated by the
pair $(E_1,E_2)$ has type $(+)$, $h=\vert\chi(E_i,E_{i+1})\vert >2$, and
for any $x\in S$ the Kronecker module $ev_{E_1,E_2,x}$ is stable.
Suppose a Mukai vector $v=(r,c,s)$
with $r>0$ satisfies the conditions $(0)$ and $(1)$ of the previous
Theorem and

\vspace{2ex}
\ldf{$(2+)$},{$\left\{ \begin{array}{rcl}\vspace{1ex}
        \mu(v)&<&\mu_{-\infty}\\
        \frac{\textstyle r(E_3)}{\textstyle r(E_2)}&<&
                             \frac{\textstyle \chi(E_3,v)}{\textstyle
\chi(E_2,v)}
\end{array} \right.
\qquad\mbox{ or }\qquad
\left\{ \begin{array}{rcl}\vspace{1ex}
      \mu_{+\infty}&<&\mu(v)\\
      \frac{\textstyle \chi(E_3,v)}{\textstyle \chi(E_2,v)}&<&
                              \frac{\textstyle r(E_1)}{\textstyle r(E_0)}
\end{array} \right.$}

\vspace{2ex}
Then for the same  $m$ and $n$ as in the previous Theorem
the Kronecker moduli space $N(h,m,n)${ }$(N^s(h,m,n))$ is isomorphic to
the coarse moduli space of semistable (resp. stable) sheaves with the
Mukai vector $v$.
\EF
\No {Proofs}
\SubNopt{Properties of the system $\{ E_i\}$}
We recall that $h=\vert\chi(E_i,E_{i+1})\vert$ does not depend on $i$ (see
2.5.1).

For $h>2$ we denote by $x_h$ the smallest solution to the equation
\begin{equation}\label{xh}
x^2-hx+1=0
\end{equation}
It is clear that $x_h<1$ and the other solution to this equation is $x_h^{-1}$.
%
%
\Pr If $h>1$, then one of the following statements hold:
\begin{description}
\item{$(+)$}
all pairs $(E_i,E_{i+1})$ are of type {\it hom\/}, all
$E_i$ are locally free;
\item{$(-)$}
there is only one {\it ext} - pair $(E_p,E_{p+1})$ , the pairs
$(E_i,E_{i+1}),i\not=p$, are {\it hom\/}, and
all $E_i,i\not=p$ are locally free.
\end{description}
\EF
\proof For any exceptional pair put by definition
$$\sgn_{E,F}=\sgn\chi_-(E,F)=\left\{ \begin{array}{rl}\vspace{.5ex}
        1&\mbox{  whenever }(E,F)\mbox{  has type {\it hom}}\\ \vspace{.5ex}
       -1&\mbox{  whenever }(E,F)\mbox{  has type {\it ext}}\\ \vspace{.5ex}
        0&\mbox{  whenever }(E,F)\mbox{  has type {\it zero}}
\end{array} \right.$$

It can be easily checked by direct calculation that
$$\sgn_{L_EF,E}\,\sgn_{E,F}\,r(L_EF)-\sgn_{E,F}\,r(E)\vert\chi(E,F)\vert
  +r(F)=0$$
(see the triples (\ref{lr}), (\ref{ls}), and (\ref{le})).
Similarly, we have analogous identity for the right mutation:
$$\sgn_{F,R_F E}\,\sgn_{E,F}\,r(R_F E)-\sgn_{E,F}\,r(F)\vert\chi(E,F)\vert
  +r(E)=0$$

Consider the system $\{ E_i\}$ generated by the pair $(E,F)=(E_1,E_2)$ and set
$$\begin{array}{llll}\vspace{1ex}
   s_1=1\ ,&s_i=\prod\limits_{1\leq k<i}\sgn_{E_k,E_{k+1}}&\mbox{for
}i>1&\mbox{and} \\
           &s_i=\prod\limits_{i\leq k<1}\sgn_{E_k,E_{k+1}}&\mbox{for }i<1&
  \end{array}$$
 From the above formulas it follows that the sequence
$$r_i=s_ir(E_i)$$
satisfies the recursion
$$r_{i-1}-hr_i+r_{i+1}=0$$
Therefore, we have two cases:
$$\begin{array}{cl}
                r_i=a+bi\qquad&\mbox{whenever }h=2\ ,\\
    r_i=Ax_h^i+Bx_h^{-i}\qquad&\mbox{whenever }h>2\ ,
\end{array}$$

\noindent\hbox to\textwidth
{where\hfil $\quad A=\frac{\textstyle{x_hr_1-r_0}}{\textstyle{x_h^2-1}}\quad ,
        \hfil B=\frac{\textstyle{x_h(x_hr_0-r_1)}}{\textstyle{x_h^2-1}} $\hfil}

\vspace{1ex}

We see that in all cases the sequence $\{ r_i\}$ has at most one change of
sign. Thus, the system $\{ E_i\}$ has at most one {\it ext} - pair.

It is readily seen that all $r_i=r(E_i)>0$ whenever $b=0$ or $AB>0$,
in this case $s_i=1,\forall i$. Therefore, all pairs in $\{ E_i\}$ are
{\it hom}, and  2.4 implies that  all $E_i$ are locally free.
In the other case, $b\not=0$ or $AB<0$, the sequence $\{ r_i\}$ is strictly
monotone. Hence, only one $r_i$ can be equal to 0. Suppose $r_p=0$;
then all $E_i,i\not=p$ are locally free and $E_p\cong {\cal O}_e(d)$.
It is clear that the unique {\it ext} - pair is just $(E_p,E_{p+1})$,
because $\hom(E_p,E_{p+1})=0$. This completes the proof.

\vspace{1ex}

Having a {\it hom} - pair $(E_i,E_{i+1})$, one can easily define the type
of the system $\{ E_i\}$ (if $(E_i,E_{i+1})$ is {\it ext\/}, then obviously,
$\{ E_i\}$ is $(-)$).

If $h=2$, then from the proof of Proposition 3.1.1 it follows that
$$r(E_i)=r(E_{i+1})\Longleftrightarrow\{ E_i\}\mbox{ is }(+)\quad\mbox{and}
   \quad r(E_i)\not=r(E_{i+1})\Longleftrightarrow\{ E_i\}\mbox{ is }(-)\ .
$$
The following statement helps to define the type of $\{ E_i\}$ provided
$h>2$.
\Pr Let $(E_i,E_{i+1})$ be a {\it hom} - pair with $h>2$; then
$$\begin{array}{cccc}
r_i^2+r_{i+1}^2-hr_ir_{i+1}<0&\Longleftrightarrow&\{ E_i\}\mbox{ is
}(+)&\mbox{and}\\
r_i^2+r_{i+1}^2-hr_ir_{i+1}>0&\Longleftrightarrow&\{ E_i\}\mbox{ is }(-)&
  \end{array}$$
\EF
\proof Without loss of generality we can assume that $i=0$.
 From the proof of Proposition 3.3.2 it follows that $\{ E_i\}$ is $(+)$ if
and only if there is no change of sign in the sequence $\{ r_i\}$, i.e.
$AB>0$. By direct calculation we obtain that
$$AB=\frac{x_h[(x_h^2+1)r_0r_1-x_h(r_0^2+r_1^2)]}{(x_h^2-1)^2}=
      \frac{x_h^2(hr_0r_1-r_0^2-r_1^2)}{(x_h^2-1)^2}
$$
This completes the proof.

\Pr Suppose $h>2$ and Mukai vector $v$ satisfies the condition $(0)$
of Theorem 2.7.1; then $\dim N(h,m,n)>0$ iff either
$$\begin{array}{ccll}\vspace{1ex}
  (2-)&\mu_{+\infty}<\mu(v)<\mu_{-\infty}&
               \mbox{whenever }\{ E_i\}\mbox{ is }(-)&\quad\mbox{or}\\
  \left( \begin{array}{c}
           \mbox{part}\\ \mbox{of }2+ \end{array} \right)&
            \begin{array}{l} \mu(v)<\mu_{-\infty}\mbox{ or }\\
                              \mu_{+\infty}<\mu(v)
              \end{array}&
 \mbox{whenever }\{ E_i\}\mbox{ is }(+)&
\end{array} $$
\EF
\proof Clearly, we can replace the rank function in the reasoning of the
proof of 3.3.2 by any additive function, in particular, to $c_1\cdot (-K)$ .
Whence, we obtain that the sequence
$$d_i=s_i\,c_1(E_i)\cdot (-K)\quad\mbox{has the form}\quad
                 d_i=Cx_h^i+Dx_h^{-i}\ ,
$$
\noindent\hbox to\textwidth
{where\hfil $\quad C=\frac{\textstyle{x_hd_1-d_0}}{\textstyle{x_h^2-1}}\quad ,
        \hfil D=\frac{\textstyle{x_h(x_hd_0-d_1)}}{\textstyle{x_h^2-1}} $\hfil
.}

\vspace{1ex}
Obviously, $\mu(E_i)=\frac{\textstyle{d_i}}{\textstyle{r_i}}$. Since $0<x_h<1$,
it follows that the sequence $\{\mu(E_i)\}$ has limits
$$\mu_{-\infty}=\frac{C}{A}=\frac{x_hd_1-d_0}{x_hr_1-r_0}\quad\mbox{and}\quad
   \mu_{+\infty}=\frac{D}{B}=\frac{x_hd_0-d_1}{x_hr_0-r_1}\ .
$$
After shift of the enumeration in $\{ E_i\}$, we obtain the formulas:
\begin{equation}\label{mu}
\mu_{-\infty}=\frac{x_hd_{i+1}-d_i}{x_hr_{i+1}-r_i}\quad\mbox{and}\quad
   \mu_{+\infty}=\frac{x_hd_i-d_{i+1}}{x_hr_i-r_{i+1}}\ ,\ \forall i
\end{equation}
Now let us consider the vector $v$ and the collection
$(E_1,E_2,F_2,\cdots,F_l)$. Since this collection is exceptional and full,
it follows that the vectors $v(E_1)$, $v(E_2)$,
$v(F_2),\cdots,v(F_l)$ form a basis in Mu$(S)$. Hence, $v$ has an
expression
\begin{equation}\label{v}
v=m'v(E_1)+n'v(E_2)+\beta_2\,v(F_2)+\cdots+\beta_l\,v(F_l)
\end{equation}
Taking the value of $\chi(F_l,\ \cdot\ )$ at the both sides,
we have that
$$0=\chi(F_l,v)=\beta_l\,\chi(F_l,F_l)=\beta_l
$$
Similarly, $\beta_j=0,\ j=l-1,\cdots ,2$. Taking the values of functions
$\chi(E_3,\ \cdot\ )$ and $\chi(E_2,\ \cdot\ )$
at the both sides of (\ref{v}), we have:
$$m'=\chi(E_3,E_1)\,\chi(E_3,v)\ ,\quad n'=\chi(E_2,v)
$$
Finally, ve obtain that
$$\mu(v)=\frac{(m'c_1(E_1)+n'c_1(E_2))\cdot (-K)}{m'r(E_1)+n'r(E_2)}\quad .
$$
Consider the linear-fractional function
$$f(x)=\frac{xd_1-d_2}{xr_1-r_2}\quad .
$$
Clearly, $\mu_{+\infty}=f(x_h),\ \mu_{-\infty}=f(x_h^{-1})$, and
$\mu(v)=f\left(\mp\frac{\textstyle{m'}}{\textstyle{n'}}\right)$ (the sign $-/+$
in
front of $\frac{\textstyle{m'}}{\textstyle{n'}}$ correspond to the type
{\it hom}/{\it ext} of the pair $(E_1,E_2)$).

\vspace{1ex}

Since $\quad f'(x)=\frac{\textstyle{d_2r_1-d_1r_2}}{\textstyle{(xr_1-r_2)^2}}=
     \frac{\textstyle{\vert\chi_-(E_1,E_2)\vert}}{\textstyle{(xr_1-r_2)^2}}>0\
,$

\vspace{1ex}

\noindent it follows that $f(x)$ is increasing for
$x>\frac{\textstyle{r_2}}{\textstyle{r_1}}$
and for $x<\frac{\textstyle{r_2}}{\textstyle{r_1}}$ .

Suppose $\{ E_i\}$ is $(+)$; then Proposition 3.3.3 implies that
$$x_h<\frac{r_2}{r_1}=\frac{r(E_2)}{r(E_1)}<x_h^{-1}\ .
$$
Thus, in this case we have:
$$\begin{array}{cccr}\vspace{1ex}
    \mu(v)<\mu_{-\infty}&\Longleftrightarrow&
      x_h<-\frac{\textstyle{m'}}{\textstyle{n'}}<
             \frac{\textstyle{r_2}}{\textstyle{r_1}}&\mbox{ and}\\
    \mu_{+\infty}<\mu(v)&\Longleftrightarrow&
      \frac{\textstyle{r_2}}{\textstyle{r_1}}<
          -\frac{\textstyle{m'}}{\textstyle{n'}}<x_h^{-1}& .
\end{array}
$$
Similarly, in the other case, when $\{ E_i\}$ is $(-)$, we obtain:
$$\mu_{+\infty}<\mu(v)<\mu_{-\infty}\quad\Longleftrightarrow \quad
                                 x_h<-\frac{m'}{n'}<x_h^{-1}\ .
$$
It can be easily proved that $\vert\chi(E_3,E_1)\vert =1$ (by applying
of the functor $\hom(F,\ \cdot\ )$ to the triples (\ref{lr}), (\ref{ls}),
and (\ref{le})). Therefore, $\vert m'\vert =m$ and $\vert n'\vert =n$.

Finally, we obtain that the conditions $(2-)$ and (part of 2+)
are equivalent to
$$x_h<\frac{m}{n}<x_h^{-1}\ .
$$
Clearly, it is equivalent to the condition
$$\dim N(h,m,n)=hnm-m^2-n^2+1>0
$$
This concludes the proof.
%
%
\SubNopt{Spectral sequence of Beilinson type}
Let $(F_0,F_1,\cdots ,F_l)$ be an exceptional collection of sheaves on $S$.
It is known (see [Go2]) that replacing the pair $(F_j,F_{j+1})$ in this
collection by one
 of the pairs $(L_{F_j}F_{j+1},F_j)$ or $(F_{j+1},R_{F_{j+1}}F_j)$ gives
always a new exceptional collection. By iterating such procedure one can
construct the {\it left dual} collection $({^\vee F}_{-l},{^\vee F}_{-l+1},
\cdots ,{^\vee F}_{-1},{^\vee F}_0)$, where
$${^\vee F}_0=F_0\quad{\rm and}\quad {^\vee F}_p=L_{F_0}\cdots
L_{F_{-p-1}}F_{-p}\quad ,
\quad p=-1,\cdots ,-l$$
\Pr Suppose $(F_0,F_1,\cdots ,F_l)$ is full exceptional collection on
Del Pezzo surface $S$. Then for any coherent sheaf $V$ on $S$ there exists
a spectral sequence with
$$E^{p,q}_1=\ext^{q-{\Delta }_p}(F_{-p},V)\otimes {^\vee F}_p\quad ,\quad
p=-l,\cdots ,0$$
where ${\Delta }_p$ is the number of non-regular mutations in obtaining
${^\vee F}_p$ from $F_{-p}$. This sequence has $E^{p,q}_\infty =0$
whenever $p+q\not= 0$ and converges to $V$ on the diagonal $p+q=0$.
\EF
This result was proved by Gorodentsev (see [Go2]).
%
\Pr Let $\;(E_1\,,\;E_2\,,\;F_2\,,\;\cdots\,,\;F_l)\;$ be full exceptional
collection, where
$F_2,\cdots,F_l$ are bundles and $\vert\chi(E_1,E_2)\vert >2$.
 Suppose a Mukai vector $v$ satisfies conditions $(0)$ and $(1)$
of Theorem 2.7.1 and $V$ is a semistable sheaf with $v(V)=v$. Then
the following statements hold:

\noindent 1. For $m=\vert\chi(E_3,v)\vert$ and
 $n=\vert\chi(E_2,v)\vert$ there exists one of the exact triples:

\vspace{1ex}

\noindent\hbox to\textwidth
{$(r)$\hfil $0\longrightarrow E_1\otimes\CC^m\longrightarrow E_2\otimes\CC^n
\longrightarrow V\longrightarrow 0$\hfil}

\vspace{1ex}

\noindent\hbox to\textwidth
{$(l)$\hfil $0\longrightarrow V\longrightarrow E_1\otimes\CC^m\longrightarrow
E_2\otimes\CC^n\longrightarrow 0$\hfil}

\vspace{1ex}

\noindent\hbox to\textwidth
{$(e)$\hfil $0\longrightarrow E_2\otimes\CC^n\longrightarrow V
\longrightarrow E_1\otimes\CC^m\longrightarrow 0$\hfil}

\vspace{1ex}

\noindent If $(r)$ takes place, then $\quad\CC^m=\hom(E_3,V)\
,\quad\CC^n=\hom(E_2,V)$.

\vspace{1ex}

\noindent 2. For the slope $\mu(V)$ there is one of the possibilities

\vspace{1ex}

a) $\exists i, \mu(V)=\mu(E_i)\quad\Longleftrightarrow\quad
                V\cong kE_i=\underbrace{E_i\oplus\cdots\oplus E_i}_{k\ times}$
b) $\mu(V)\not\in\{\mu(E_i)\}$, then

\vspace{1ex}

$\mu(V)<\mu_{-\infty}$ or $\mu_{+\infty}<\mu(V)$ whenever $\{ E_i\}$ is of type
$(+)$
and

\vspace{1ex}

$\mu_{+\infty}<\mu(V)<\mu_{-\infty}\qquad$ whenever $\{ E_i\}$ is of type
$(-)$.

\vspace{1ex}

\noindent 3. If $\mu(V)\not\in\{\mu(E_i)\}$ then the existence of the triples
$(r)$, $(l)$, and $(e)$ is characterized as follows:

\vspace{1.7ex}
\noindent\hbox to\textwidth
{\begin{tabular}{lcl}
\vspace{1.7ex}
$(r)$&$\Longleftrightarrow$&
          \begin{tabular}{l}
            $\{E_i\}$ has type $(+)$ and $\mu_{+\infty}<\mu(V)$\\
              or $\{E_i\}$ has type $(-)$ and the {\it ext} - pair
             $(E_p,E_{p+1})$ has $p<1$
           \end{tabular}\hfil\\

\vspace{1.7ex}
$(l)$&$\Longleftrightarrow$&
          \begin{tabular}{l}
            $\{E_i\}$ has type $(+)$ and $\mu(V)<\mu_{-\infty}$\\
              or $\{E_i\}$ has type $(-)$ and the {\it ext} - pair
             $(E_p,E_{p+1})$ has $p>1$
           \end{tabular}\hfil\\

$(e)$&$\Longleftrightarrow$& $\ \{E_i\}$ has type $(-)$ and the {\it ext} -
pair is $(E_1,E_2)$.
\end{tabular}}
\EF
%
\proof
The condition $(1)$ implies that $\mu(F_j\otimes\omega_S)<\mu(V)<\mu(F_j),
\forall j$. From the stability properties and Serre
duality it follows that $\hom (F_j,V)=0$ and
$\ext^2(F_j,V)\cong\hom^*(V,F_j\otimes\omega_S)=0$. Using $(0)$, we have
$\ext^1(F_j,V)=0$. Therefore,  $\ext^i(F_j,V)=0,\forall i,j$.

Consider the spectral sequence associated with the full exceptional
collection
 $(E_2,E_3,F_2,\cdots,F_l)$. From the above it follows that
$E^{p,q}_1=0$ for $p\leq -2$. Suppose for definiteness that the left mutation
of
$(E_2,E_3)$ is regular ($\Leftrightarrow\ $ the right mutation of
$(E_1,E_2)$ is regular).
Then $(E_1,E_2)$ is a {\it hom} - pair and
the $E_1$-term of the spectral sequence has the form
$$\begin{array}{cccccc}
      & & & & \llap{$\scriptstyle q$}\bigl\uparrow\bigr. & \\
     0&\cdots &0&0&0& \\
     \cdots &\cdots &\cdots &\cdots &\cdots & \\ \vspace{1ex}
     0&\cdots &0&0&0& \\ \vspace{1ex}
     0&\cdots &0&E_1\otimes\ext^2(E_3,V)&E_2\otimes\ext^2(E_2,V)& \\
\vspace{1ex}
     0&\cdots &0&E_1\otimes\ext^1(E_3,V)&E_2\otimes\ext^1(E_2,V)& \\
     0&\cdots &0&E_1\otimes\hom (E_3,V) &E_2\otimes\hom (E_2,V) &
                                     \stackrel{p}{\longrightarrow} \\
    \end{array}$$

It is readily seen that $E^{p,q}_1=E^{p,q}_\infty =0$ for $q=2$, $p=0,-1$.
Whence we have $\ext^2(E_3,V)=\ext^2(E_2,V)=0$.

The convergence to $V$ as $p+q=0$ implies that there are the
following exact triples:
$$0\longrightarrow E_1\otimes\hom (E_3,V)\longrightarrow E_2\otimes\hom (E_2,V)
\longrightarrow {\cal F}_0\longrightarrow 0
$$
$$0\longrightarrow {\cal F}_0\longrightarrow V
\longrightarrow {\cal F}_1\longrightarrow 0
$$
$$0\longrightarrow {\cal F}_1 \longrightarrow
E_1\otimes\ext^1(E_3,V)\longrightarrow
E_2\otimes\ext^1(E_2,V)\longrightarrow 0
$$

If $\mu(V)>\mu(E_1)$, then from the stability properties it follows that
$\hom(V,E_1)=0$. Therefore, ${\cal F}_1=0$, ${\cal F}_0=V$ and
$\ext^1(E_3,V)=\ext^1(E_2,V)=0$. Finally, we have the triple $(r)$, and
from the Swing Lemma it follows that $\mu(V)>\mu(E_2)$.

 From the other hand, if $\mu(V)\leq\mu(E_1)<\mu(E_2)$, then $\hom(E_2,V)=0$,
and we obtain: ${\cal F}_0=0$, ${\cal F}_1=V$, and $\hom (E_3,V)=\hom
(E_2,V)=0$.

It is obvious that in the both cases the dimensions of the vector
spaces that are tensored with the bundles $E_1$ and $E_2$  equal
respectively $\vert\chi(E_3,v)\vert$ and $\vert\chi(E_2,v)\vert$.

\vspace{1ex}

If the right mutation of $(E_1,E_2)$ is singular (or extension), then, applying
the spectral sequence associated with the collection
$(E_2,E_3,F_2,\cdots,F_l)$, one
can easily check in the similar way that $(l)$ (respectively, $(e)$) takes
place. This completes the proof of 1.

\vspace{1ex}

Now we prove 2.

\noindent a) One can apply 1. to the full collection
$(E_i,E_{i+1},F_2,\cdots,F_l)$, so,
without loss of generality we can renumerate $\{ E_i\}$ such that
$\mu(V)=\mu(E_1)$. Thus we have one of the triples  $(r)$, $(l)$, or $(e)$,
and from the Swing Lemma it immediately follows that $V\cong E_1\otimes\CC^m$.

\noindent b) Actually, we have proved the following lemma.
%
%
\Lm If the right mutation of $(E_i,E_{i+1})$ is regular, then
$$\mu(E_i)<\mu(V)\ \Rightarrow\ \mu(E_{i+1})<\mu(V)\quad{\sl and }\quad
\mu(V)<\mu(E_{i+1})\ \Rightarrow\ \mu(V)<\mu(E_i)\ .
$$
\EF
Suppose $\{ E_i\}$ is of type $(+)$; then all mutations of pairs are regular
and we obviously have one of the possibilities: $\mu(E_i)<\mu(V),\forall i$ or
$\mu(V)<\mu(E_i),\forall i$. Consequently, $\mu_{+\infty}\leq\mu(V)$ or
$\mu(V)\leq\mu_{-\infty}$. From the explicit formulas (\ref{mu}) it is
not hard to check that $\mu_{+\infty}$ and $\mu_{-\infty}$ are irrational,
so the inequalities are strict.

Suppose $\{ E_i\}$ is of type $(-)$ and the {\it ext} - pair is
$(E_p,E_{p+1})$.
According to the proof of 1., we have the exact triples:
$$0\longrightarrow V\longrightarrow E_{p-1}\otimes\CC^{m_{p-1}}\longrightarrow
E_p\otimes\CC^{n_{p-1}}\longrightarrow 0
$$
$${\rm and}\quad 0\longrightarrow E_{p+1}\otimes\CC^{n_p}\longrightarrow V
\longrightarrow E_p\otimes\CC^{m_p}\longrightarrow 0
$$
Therefore, $\mu(E_{p+1})<\mu(v)<\mu(E_{p-1})$. Applying Lemma 3.2.3,
we have: $\mu(E_i)<\mu(V)$ for $i\geq p+1$ and $\mu(V)<\mu(E_i)$ for $i\leq
p-1$.
As in the previous case, using 2.5, we obtain:
$\mu_{+\infty}<\mu(V)<\mu_{-\infty}$.

The statement 3. immediately follows from 2.b) and the Swing Lemma.
This concludes the proof of Proposition 3.2.2.
%
%
\SubNopt{Correspondence between sheaves and Kronecker modules}

Further on through this section we shall assume that
$(E_1,E_2,F_2,\cdots ,F_l)$ is a full exceptional collection,
$\{ E_i\} $ is the system of sheaves generated by $(E_1,E_2)$,
$v$ is a Mukai vector satisfying conditions $(0)$ and $(1)$ of
Theorem 2.7.1, and $V$ is a semistable sheaf with $v(V)=v$.

Proposition 3.2.2 guarantees the existence of one of the triples
$(r),\ (l)$ or $(e)$. For definiteness we suppose that the triple

\vspace{2ex}

\noindent\hbox to\textwidth
{$(r)$\hfil $0\longrightarrow E_1\otimes\CC^m\stackrel{a}{\longrightarrow}
E_2\otimes\CC^n\longrightarrow V\longrightarrow 0$\hfil}

\vspace{2ex}

\noindent takes place (the case of existing $(l)$ is dual)\footnote
{\ If there exists $(e)$ then $(E_1,E_2)$ is an {\it ext} - pair and
one can shift the numeration in the system $\{ E_i\}$ such that the
new $(E_1,E_2)$ will be {\it hom}.}.
Here $a$ is the differential $d_1$ of the spectral sequence associated
with the collection $(E_2,E_3,F_2,\cdots ,F_l)$,
$$\CC^m=\hom(E_3,V)\ ,\quad\CC^n=\hom(E_2,V)
$$
(see proof of 3.2.2). It is clear that $\hom(E_1,E_2)\not= 0$ and the pair
$(E_1,E_2)$ is of type {\it hom} (see 2.4). Hence, we have
$$\chi_-(E_1,E_2)>0\quad {\rm and}\quad\mu (E_1)<\mu(E_2)
$$
(from the statement 3. of Proposition 3.2.2 and 2.5 it follows that
$E_1$ and $E_2$ are locally free).

\vspace{1ex}
We assign to $V$ the Kronecker $\hom^*(E_1,E_2)$ - module
\begin{equation}\label{tv}\MAP{t=t_V},\CC^m\otimes\hom^*{(E_1,E_2)},\CC^n
\end{equation}
\noindent naturally associated with the morphism of bundles $a$ in $(r)$.
The aim of this subsection is to show that the sheaf $V$ is (semi)stable if and
only if $t$ is (semi)stable as Kronecker module (see 2.6).
\Pr \\
a) $t$ is not semistable $\quad\Longrightarrow\quad V$ is not semistable;\\
b) $t$ is semistable but not stable $\quad\Longrightarrow\quad V$ is not
stable.
\EF
\proof
Let a pair $(f_0,f_1)$ of linear maps, $\MAP {f_0},{H'_0},{\CC^m}$ ,
 $\MAP {f_1},{H'_1},{\CC^n}$ , define a Kronecker submodule
$\MAP {t'},{H'_0\otimes\hom^*{(E_1,E_2)}},{H'_1}$ of $t$ and let $t''=t/t'$
be the Kronecker quotient-module.
Consider the following commutative diagram:
$$\begin{array}{ccccccccc}
0&\longrightarrow &E_1\otimes H'_0&\stackrel{{\rm id}_{E_1}\otimes f_0}
 {\longrightarrow}&E_1\otimes \CC^m&\longrightarrow&E_1\otimes(\CC^m/H'_0)&
  \longrightarrow 0  \\
 & &\Biggl\downarrow\Biggr.\rlap{$\scriptstyle a'$}& &
     \Biggl\downarrow\Biggr.\rlap{$\scriptstyle a$}& &
      \Biggl\downarrow\Biggr.\rlap{$\scriptstyle a''$}& & \\
0&\longrightarrow &E_2\otimes H'_1&\stackrel{{\rm id}_{E_2}\otimes f_1}
 {\longrightarrow}&E_2\otimes \CC^n&\longrightarrow&E_2\otimes(\CC^n/H'_1)&
  \longrightarrow 0  \\
\end{array},$$
where $a'$ and $a''$ naturally correspond to $t'$ and $t''$. The morphism
$a'$ is injective as restriction of $a$. Hence, if we regard this
diagram as an exact triple of complexers, then the long cohomology sequence
for it has the form:
$$0\longrightarrow{\cal K}\stackrel{\varphi}{\longrightarrow}V'
    \stackrel{\psi}{\longrightarrow}V\longrightarrow{\cal C}
     \longrightarrow 0\ ,
$$
where ${\cal K}=ker\ a'',\ V'=coker\ a',\ {\cal C}=coker\ a''$. Let us
note that $V'\not= 0$ because $E_1\otimes H'_0\not\cong E_2\otimes H'_1$.

It can be easily seen that
$$\left| \begin{array}{cc} \vspace{1ex}
           c_1(V')\cdot (-K)& r(V') \\
           c_1(V)\cdot (-K)& r(V)
  \end{array} \right| =
\left| \begin{array}{cc} \vspace{1ex}
           \dim H'_1 & \dim H'_0 \\
                 n   &    m
  \end{array} \right| \cdot
\left| \begin{array}{cc} \vspace{1ex}
           c_1(E_2)\cdot (-K)& r(E_2) \\
           -c_1(E_1)\cdot -(-K)& r(E_1)
  \end{array} \right|
$$
This implies that
$$\chi_-(V,V')=-\left| \begin{array}{cc} \vspace{1ex}
           \dim H'_1 & \dim H'_0 \\
                 n   &    m
  \end{array} \right| \cdot
\chi_-(E_1,E_2)\quad .
$$
Further reasonings for the cases a) and b) are parallel.
Suppose $t'$ violates the semistability (resp. stability) criterion
2.6.1. That means that the determinant in the last equality is $<0\
(\leq 0)$. Therefore, we have
$$\chi_-(V,V')>0\ (\geq 0)\quad .$$

If ${\cal K}\not= 0$ then $r({\cal K})\not= 0$ and $\mu({\cal K})\leq\mu(E_1)$
because {\cal K} is a subsheaf of the semistable sheaf
$E_1\otimes(\CC^m/H'_0)$.
Hence, we have $r(v')\not= 0$ and
$$\mu(V')>(\geq )\ \mu(V)>\mu(E_1)\geq\mu({\cal K})\quad .
$$
 From the Swing Lemma it follows that $\mu(U)>(\geq )\ \mu(V')>(\geq )\
\mu(V)$,
where $U=coker\  \psi = ker\  \varphi$ (since $U$ is a non-zero subsheaf of a
torsion free sheaf $V$, we see that $r(U)\not= 0$).
If ${\cal K}=0$ then we obviously have $U\cong V'$. Finally we obtain that
in the case a) $V$ is not semistable and in the case b) $V$ is not stable.
This completes the proof.
\Pr \\
a) $V$ is not semistable $\quad\Longrightarrow\quad t$ is not semistable;\\
b) $V$ is semistable but not stable $\quad\Longrightarrow\quad t$ is not
stable.
\EF
First we prove the following statement.
\Lm There exists a subsheaf $U\subset V$ with $\mu(U)>\mu(V)$ in the case a)
and $\mu(U)\geq\mu(V)$ in the case b) satisfying the properties:

\vspace{1ex}

\noindent{\rm (i)} if $E$ is a stable sheaf with $\chi(E,V)=0$ and
$\mu(E)-K^2<\mu(V)<\mu(E)$ , then
$$\ext^i(E,U)=0\ ,\ \forall i\ ;$$
\noindent{\rm (ii)} if $F$ is a stable sheaf and $\mu(F)<\mu(V)$, then
$\quad\hom(U,F)=0$.
\EF
\proof a) Suppose $V$ is not semistable. Consider the Harder -- Narasimhan
filtration of $V$
$$0={\cal F}_{n+1}\subset{\cal F}_n\subset\cdots\subset{\cal F}_2\subset{\cal
F}_1=V
$$
This is a unique filtration such that all the quotients $G_i={\cal F}_i/{\cal
F}_{i+1}$ are semistable, $\mu(G_{i+1})>\mu(G_i)$ for $i=1,\cdots,n-1$, and
$\mu({\cal F}_i)<\mu({\cal F}_{i+1}),\ i=1,\cdots,n$.\footnote{\ This
filtration
is constructed by induction: $G_i$ is the torsion free quotient-sheaf of
${\cal F}_i$ with the smallest slope and maximal rank (see, for example,
[Ku]).}
Since $\mu(G_n)=\mu({\cal F}_n)>\mu(V)>\mu(V/{\cal F}_2)=\mu(G_1)$ (the last
inequality holds by the Swing Lemma), we conclude that there exists a unique
$p$ such that

\vspace{2ex}

\rdf{($\ast$)},{$\mu(G_1)<\cdots <\mu(G_{p-1})\leq\mu(V)\leq\mu(G_p)
      <\cdots <\mu(G_n)$}

\vspace{2ex}

We put $U={\cal F}_p$. The inequality $\mu(U)>\mu(V)$ holds by the properties
of the filtration.

Now we prove (i).
Let us note that $\ext^i(E,V)=0,\forall i$. Really, from the Lemma conditions
and stability properties it follows that $\hom(E,V)=0$ and $\ext^2(E,V)\cong
\hom^*(E,V)=0$. Using the condition $\chi(E,V)=0$, we obtain that
$\ext^1(E,V)=0$.

Consider the exact triples

\vspace{2ex}

\rdf{($\ast\ast$)},{$
0\longrightarrow {\cal F}_{i+1} \longrightarrow {\cal F}_i\longrightarrow G_i
\longrightarrow 0$}

\vspace{2ex}

We consequently apply the functor $\hom(E,\ \cdot\ )$ to these triples starting
from
 $i=1$ to $i=p-1$. From the inequality ($\ast$), stability of $E$ and $G_i$,
and the Lemma conditions it follows that $\hom(E,G_i)=0,\ i=1,\cdots,p-1$.
Whence, on each step we have the implication
$$\hom(E,{\cal F}_i)=\ext^1(E,{\cal F}_i)=0\quad\Longrightarrow\quad
  \hom(E,{\cal F}_{i+1})=\ext^1(E,{\cal F}_{i+1})=0
$$
\noindent (for $i=1$ the former equality holds because ${\cal F}_1=V$). On the
last step we obtain: $\hom(E,{\cal F}_p)=\ext^1(E,{\cal F}_p)=0$.

Similarly, from ($\ast$), stability of $E\otimes K$ and $G_i$,
and the Lemma conditions it follows that $\hom(G_i,E\otimes K)=0,\
i=p,\cdots,n$.
Applying consequently the functor $\hom(\ \cdot\ ,E\otimes K)$ to the triples
($\ast\ast$) starting from $i=n$ to $i=p$, we have on each step the implication
$$\hom({\cal F}_{i+1},E\otimes K)=0\quad\Longrightarrow\quad
  \hom({\cal F}_i,E\otimes K)=0
$$
\noindent Thus, for $i=p$ we obtain: $\hom^*({\cal F}_p,E\otimes K)=
\ext^2(E,{\cal F}_p)=0$.

In order to prove (ii) we consequently apply the functor $\hom(\ \cdot\ ,F)$
to the triples ($\ast\ast$) starting from $i=n$ to $i=p$. As above, on the
last step we get $\hom(U,F)=0$.

\noindent b) Suppose $V$ is semistable but not stable; then by definition
there is a subsheaf $U\subset V$ such that $0<r(U)<r(V)$ and $\mu(U)=\mu(V)$.
It is clear that $U$ is semistable. Whence, (i) and (ii) follow from the
stability properties.

This completes the proof of Lemma.
\par\noindent{\sc Proof of Proposition 3.3.2.}\hspace{1ex}
Consider the subsheaf $U\subset V$ as in the previous Lemma. Applying the
statement
(i) to $E=F_j$, we have: $\ext^i(F_j,U)=0,\forall i,\ j=2,\cdots ,l$. Then the
spectral sequence associated with the collection $(E_2,E_3,F_2,\cdots,F_l)$
and convergent to $U$ has the $E_1$-term of the form:
$$\begin{array}{cccccc}
      & & & & \llap{$\scriptstyle q$}\bigl\uparrow\bigr. & \\
     0&\cdots &0&0&0& \\
     \cdots &\cdots &\cdots &\cdots &\cdots & \\ \vspace{1ex}
     0&\cdots &0&0&0& \\ \vspace{1ex}
     0&\cdots &0&E_1\otimes\ext^2(E_3,U)&E_2\otimes\ext^2(E_2,U)& \\
\vspace{1ex}
     0&\cdots &0&E_1\otimes\ext^1(E_3,U)&E_2\otimes\ext^1(E_2,U)& \\
     0&\cdots &0&E_1\otimes\hom (E_3,U) &E_2\otimes\hom (E_2,U) &
                                     \stackrel{p}{\longrightarrow} \\
    \end{array}$$

It is readily seen that $E^{p,q}_1=E^{p,q}_\infty =0$ for $q=2$, $p=0,-1$.
Whence we have $\ext^2(E_3,U)=\ext^2(E_2,U)=0$.

The convergence to $U$ as $p+q=0$ implies that there are the
following exact triples:
$$0\longrightarrow E_1\otimes\hom (E_3,U)\longrightarrow E_2\otimes\hom (E_2,U)
\longrightarrow U_0\longrightarrow 0
$$
$$0\longrightarrow U_0\longrightarrow U
\longrightarrow U_1\longrightarrow 0
$$
$$0\longrightarrow U_1 \longrightarrow E_1\otimes\ext^1(E_3,U)\longrightarrow
E_2\otimes\ext^1(E_2,U)\longrightarrow 0
$$
Let us recall that $\mu(E_1)<\mu(V)$. By the statement (ii) of the previous
Lemma
we have: $\hom(U,E_1)=0$. From the last exact triples it follows that
$U_1=0$, $\ext^1(E_3,U)=\ext^1(E_2,U)=0$, and $U\cong U_0$.

We see that the subsheaf $U$ can be included in the exact triple
$$0\longrightarrow E_1\otimes H'_0\longrightarrow E_2\otimes
H'_1\longrightarrow
   U\longrightarrow 0\ ,\mbox{ where}
$$
$$H'_0=\hom(E_3,U)\ ,\quad H'_1=\hom(E_2,U)
$$
Applying the functors $\hom(E_3,\ \cdot\ )$ and $\hom(E_2,\ \cdot\ )$ to
the exact triple
$$0\longrightarrow U\longrightarrow V\longrightarrow
  V/U\longrightarrow 0
$$
we have the exact sequences of vector spaces
$$0\longrightarrow H'_0\stackrel{f_0}{\longrightarrow}\CC^m=\hom(E_3,V)
\longrightarrow\cdots
$$
$$0\longrightarrow H'_1\stackrel{f_1}{\longrightarrow}\CC^n=\hom(E_2,V)
\longrightarrow\cdots
$$
The pair of embeddings $(f_0,f_1)$ defines a Kronecker submodule
$\MAP {t'},{H'_0\otimes L},{H'_1}$ of the module $t$.

As in the proof of Proposition 3.3.1, we have the formula
$$\chi_-(V,U)=-\left| \begin{array}{cc} \vspace{1ex}
           \dim H'_1 & \dim H'_0 \\
                 n   &    m
  \end{array} \right| \cdot
\chi_-(E_1,E_2)
$$
Since $\chi_-(V,U)>0(=0)$ in the case a) (respectively, b)) and
$\chi_-(E_1,E_2)>0$; we obtain that
$$\left| \begin{array}{cc} \vspace{1ex}
           \dim H'_1 & \dim H'_0 \\
                 n   &    m
  \end{array} \right| <0\ (=0)
$$
Thus, the submodule $t'$ violates the semistability (resp. stability)
criterion 2.6.1. This concludes the proof of Proposition 3.3.2.
%
%
\SubNo{} Now let us consider the full exceptional collection $(E_1,E_2,F_2,
\cdots F_l)$, the system $\{ E_i\}$ generated by $(E_1,E_2)$ such that all
conditions of one of Theorems 2.7.1 or 2.7.2 are satisfied. We recall
that $h=\vert\chi(E_1,E_2)\vert >2,\ m=\vert\chi(E_3,v)\vert ,\
 n=\vert\chi(E_2,v)\vert$. Suppose $(E_1,E_2)$ is of type {\it hom\/},
$$\MAP t,\CC^m\otimes\hom^*{(E_1,E_2)},\CC^n
$$
is Kronecker $\hom^*(E_1,E_2)$-module, and
$$\MAP a,E_1\otimes\CC^m,E_2\otimes\CC^n
$$
is the corresponding morpfism of sheaves.
\Pr If either $\{ E_i\}$ is $(-)$ and the {\it ext} - pair is $(E_0,E_1)$ or
$\{ E_i\}$ is $(+)$ and $\qquad\left\{ \begin{array}{rcl}\vspace{1ex}
      \mu_{+\infty}&<&\mu(v)\\
      \frac{\textstyle m}{\textstyle n}&<&
                              \frac{\textstyle r(E_1)}{\textstyle r(E_0)}
\end{array} \right.\ ,\quad$
then
$$t\mbox{ is (semi)stable}\quad\Longleftrightarrow\quad a\mbox{ is injective
and }\
V=coker\, a\ \mbox{ is (semi)stable}$$
\EF
\proof The implication
$$t\mbox{ is semistable}\quad\Longrightarrow\quad a\mbox{ is injective}$$
is proved in the work [Ka] for the case $S=\PP^1\times\PP^1$. For a Del
Pezzo surface $S$ the proof is similar, and it is omitted. Now the Proposition
follows from 3.3.1 and 3.3.2.
\newpage
\SubNo{Proof of main results} is by standard reasoning.

Without loss of generality we can assume that the {\it ext} - pair is
$(E_0,E_1)$
provided $\{ E_i\}$ is $(-)$. If $\{ E_i\}$ is $(+)$, we give the proof only
for the case
$$\left\{ \begin{array}{rcl}\vspace{1ex}
      \mu_{+\infty}&<&\mu(v)\\
      \frac{\textstyle \chi(E_3,v)}{\textstyle \chi(E_2,v)}&<&
                              \frac{\textstyle r(E_1)}{\textstyle r(E_0)}
\end{array} \right.
$$
the other case is dual.

Let $X$ be an algebraic variety, and let ${\cal U}$ be a coherent sheaf on
$X\times S$ such that ${\cal U}$ is f\/lat over $X$ and  ${\cal U}_x$ is a
semistable sheaf on $S$ with Mukai vector $v$ for any closed point
$x\in X$.

Let $p_X$ and $p_S$ denote the projections of $X\times S$ onto $X$ and $S$,
respectively; then the sheaves
$$H_0=R^0p_{X\ast}(p_S^\ast E_3^*\otimes {\cal U})\quad\mbox{and}\quad
  H_1=R^0p_{X\ast}(p_S^\ast E_2^*\otimes {\cal U})$$
are locally free, and there exists a canonical morpfism
$$\MAP\theta ,H_0\otimes\hom^*{(E_1,E_2)},H_1$$
 such that $\theta_x=t_{{\cal U}_x}$
for any closed point $x\in X$.

We consider a sufficiently small neighbourhood $X_x$ of the point $x$ and
trivializations $H_0\cong {\cal O}_{X_x}\otimes\CC^m$ and
$H_1\cong {\cal O}_{X_x}\otimes\CC^n$. Proposition 3.4.1 gives a morphism
$\MAP\lambda_x,X_x,W^{ss}$. By gluing together the morphisms $\lambda_x$
composed with the projection $\MAP\pi ,W^{ss},N({h,m,n})$, we obtain
a morphism $\MAP\phi_{\cal U},X,N({h,m,n})$, which depends functorially on
the classes ${\cal U}$.

Next suppose that $M$ is an algebraic variety such that, for any family
${\cal U}$ of sheaves on $S$ that is parametrized by the algebraic variety
$X$, there exists a morphism $\MAP\psi_{\cal U},X,M$ which depends functorially
on
the classes ${\cal U}$. We consider the universal Kronecker module over
$W^{ss}$:
$$T\,\colon\;{\cal O}_{W^{ss}}\otimes\CC^m\otimes\hom^*(E_1,E_2)\longrightarrow
                                                   {\cal
O}_{W^{ss}}\otimes\CC^n \ .
$$
According to 3.4.1, it gives a family ${\cal V}$ of semistable sheaves on $S$.
The induced morphism $\MAP\psi_{\cal V},W^{ss},M$ is $G_0$-invariant because of
the
functorial dependence $\psi_{\cal U}$ on ${\cal U}$. Hence, $\psi_{\cal V}$
induces a morphism $\MAP\eta ,N({h,m,n}),M$. We show that
 $\eta\circ\phi_{\cal U}=\psi_{\cal U}$. Because of functoriality of
$\phi_{\cal U}$ and $\psi_{\cal U}$, it is sufficient to show the following:
if a sheaf $V$ on $S$ is a cokernel of a sheaf embedding $\MAP
a,E_1\otimes\CC^m,E_2\otimes\CC^n$, then the Kronecker module $t$ which
canonically corresponds to $a$ is isomorphic to $t_V$. But it follows from
[Ka,4.4.2]. This completes the proof.

\newpage\par\noindent
\begin{flushleft}
 {\large\bf References.}
\end{flushleft}
\vspace{1ex}
\begin{description}
\item{[DrLP]:}
J.-M.Drezet, J.Le Potier
{\it Fibr\'es stables et fibr\'es exceptionnels sup $\PP^2$.\/}
Ann. Sci. \'Ecole Norm. Sup. (4) {\bf 18} (1985), 193-243.
\item{[Dr1]:}
J.-M.Drezet.
{\it Fibr\'es exceptionnels et suite spectrale de Beilinson g\'en\'eralis\'ee
sur
$\PP^2(\CC)$.\/}
Math. Ann. {\bf 275} (1986), 25-48.
\item{[Dr2]:}
J.-M.Drezet.
{\it Fibr\'es exceptionnels et vari\'et\'es de modules de faisceaux
semi-stables sur
$\PP^2(\CC)$.\/}
J. Reine Angew. Math. {\bf 380} (1987), 14-58.
\item{[Dr3]:}
J.-M.Drezet.
{\it Vari\'et\'es de modules extr\'emales de faisceaux semi-stables sur
$\PP^2(\CC)$.\/}
Math. Ann. {\bf 290} (1991), 727-770.
\item{[Go1]:}
A.L.Gorodentsev.
{\it Surgeries of exceptional vector bundles on $\PP_n$.\/}
English transl. in Math. USSR Izv. {\bf 32} (1989) No.1, 1-13.
\item{[Go2]:}
A.L.Gorodentsev.
{\it Exceptional vector bundles on surfaces with moving anticanonical
divisor.\/}
English transl. in Math. USSR Izv. {\bf 33} (1989) No.1, 67-83.
\item{[Go3]:}
A.L.Gorodentsev.
{\it Helix theory and non-symmetric bilinear forms.\/}
In: {\it Algebraic Geometry and its Applications\/} (proceedings of 8th Alg.
Geom. Conf., Yaroslavl' 1992). Aspects of Math. (1994).
\item{[Ka]:}
B.V.Karpov.
{\it Semistable sheaves on a two-dimensional quadric, and Kronecker modules.\/}
English transl. in Russian Acad. Sci. Izv. Math. Vol. {\bf 40} (1993), No.1,
33-66
\item{[KuOr]:}
S.A.Kuleshov, D.O.Orlov.
{\it Exceptional sheaves on Del Pezzo surfaces.}
Izv. Russ. Acad. Nauk Ser. Mat. Vol. {\bf 58} (1994), No.3, 59-93 (in Russian)
\item{[Ku]:}
S.A.Kuleshov.
{\it Exceptional and rigid sheaves on surfaces with anticanonical class without
base components.}
Preprint No.1 of Math. College of the Independent Moscow University, 1994 (in
Russian)
\item{[MuFo]:}
D. Mumford, J. Fogarty.
{\it Geometric invariant theory.} (Ergeb. Math. Grenzgeb. Bd.34)
Berlin Heidelberg New York: Springer 1982.

\end{description}

\end{document}